\documentclass{article}
\usepackage{pifont}
\usepackage{amsthm}
\usepackage{bbm}
\usepackage[title,titletoc]{appendix}
\usepackage{graphicx}
\usepackage{amsmath,amssymb,amsthm,amsfonts}
\usepackage{amssymb}
\usepackage{color}
\usepackage{amsfonts}
\usepackage{amssymb}
\usepackage{amsthm}
\usepackage{amsmath}
\usepackage{empheq}
\usepackage{indentfirst}
\usepackage{cite}
\usepackage{mathrsfs}
\allowdisplaybreaks[4]

\newtheorem{thm}{Theorem}[section]
\newtheorem{lem}{Lemma}[section]

\theoremstyle{definition}

\numberwithin{equation}{section}

\DeclareMathSymbol{\C}{\mathalpha}{AMSb}{"43} \topmargin-.1in

\def\om{\omega}

\def\f{\frac}
\def\F{\displaystyle\frac}
\def\dis{\displaystyle}

\def\q{\quad}
\def\ga{\gamma}
\def\I{\dis\int}
\def\i{\int}
\def\n{\nabla}
\newcommand{\beq}{\begin{equation}}
\newcommand{\eeq}{\end{equation}}
\newcommand{\bea}{\begin{array}}
\newcommand{\eea}{\end{array}}

\newcommand{\eps}{\varepsilon}

\newcommand{\lam}{\lambda}

\newcommand{\sig}{\sigma}

\def\th{\theta}
\newcommand{\pa}{\partial}
\def\dx{\mathrm{d}x}

\def\ni{\noindent}
\def\proof{{\ni\bf Proof:\quad}}
\def\pfthm2{{\ni\bf Proof of Theorem \ref{thm2}: }}
\def\pfthm3{{\ni\bf Proof of Theorem \ref{thm3}: }}
\def\proofend{{\hfill$\Box$}\\}

\title{Dynamics of a predator-prey system in open advective heterogeneous environments}
\author{Qi Wang
\thanks{College of Science, University of Shanghai for Science and Technology,
Shanghai 200093, P.R. China.
Email: \texttt{qwang@usst.edu.cn}.
}
 }
\date{}

\begin{document}

\maketitle

\begin{abstract}
In this paper, we investigate the effect of dispersal and advection on the dynamics of a predator-prey model.
More precisely,
we show that the linear stability of the semi-trivial steady state is determined by the dispersal rate, the mortality rate of the predator and the advection rate.
We point out that compared to homogeneous intrinsic growth rate and carrying capacity,
the case in this paper is more complicated.
This work gives a investigation to an open problem proposed by Nie et al. in \cite{NWW}
by considering a more general model, and then,
can be seen as a further development of their work \cite{NWW}.
\end{abstract}



\section{Introduction}

In this paper,
we study the following predator-prey system in open advective environments,
where the predator and prey are exposed to unidirectional water flow
\beq
\tag{P}
\label{P}
\left\{\arraycolsep=1.5pt
\begin{array}{lll}
u_t=d_1u_{xx}-q_1u_x+f(x,u,P)u,&x\in(0,L),t>0,\\[2mm]
P_t=d_2P_{xx}-q_2P_x+g(u)P,&x\in(0,L),t>0,\\[2mm]
d_1 u_x(0,t)-q_1u(0,t)=0,\q u_x(L,t)=0,&t>0, \\[2mm]
d_2 P_x(0,t)-q_2P(0,t)=0,\q P_x(L,t)=0,&t>0,\\[2mm]
u(x,0)=u_0(x)\ge,\not\equiv0,\q
P(x,0)=P_0(x)\ge,\not\equiv0,&x\in(0,L),
\end{array}
\right.
\eeq
where $u(x,t)$ and $P(x,t)$ are the population densities of the prey and predator species
at the time $t$ and location $x$, respectively.
Here $d_1 $ and $d_2 $ are dispersal rates of the prey and predator, respectively.
$q_1,q_2$ denote the effective advection rate caused by unidirectional water flow,
$L$ is the domain length.
$f$ and $g$ are given by
$f(x,u,P)=r(x)(1-\F{u}{K(x)})-aP, g(u)=eau-\ga$,
where $r(x)$ is the intrinsic growth rate and $K(x)$ is the carrying capacity of the prey(which reflects the resources available at
location $x\in[0,L]$).
$a$ is the predation rate, $e$ is the trophic conversion efficiency
and $\ga$ is the mortality rate of the predator.

This system was first proposed
by  Hilker and Lewis in \cite{HiLe}.
Based on model \eqref{P},
Hilker and Lewis \cite{HiLe} distinguished the flow speed scenarios,
including coexistence, survival of one species and extinction of both populations.

Recall that when $q_1=q_2=q$, $r(x)\equiv r$ and $K(x)\equiv K$ are two constants,
the system \eqref{P} has been studied in \cite{NWW},
where they considered the global stability of the trivial and semi-trivial
steady states.
More precisely, it was  showed  in \cite{NWW} that for any $d_1 ,d_2 >0$
there exist $0<q_0<q^*$ such that
$(0,0)$ is globally asymptotically stable
(among all nonnegative and nontrivial initial data)
provided that $q\ge q^*$;
while $(\vartheta,0)$ is globally asymptotically stable
if $\ga\ge br,0<q<q^*$ or $\ga<br,q_0<q<q^*$,
where $\vartheta$ is the unique positive solution of the problem
\beq
\left\{\arraycolsep=1.5pt
\begin{array}{l}
d_1  \vartheta_{xx}-q \vartheta_x+\vartheta(r-\vartheta)=0,\q x\in(0,L),\\[2mm]
d_1 \vartheta_x(0)-q \vartheta(0)=\vartheta_x(L)=0.
\end{array}\right.
\eeq
In addition,
H. Nie, B. Wang and J. H. Wu (\cite{NWW}) proposed
some open problems:
for example,
if the intrinsic growth rate or the carrying capacity of the prey are not positive constants,
if $q_1\neq q_2$,
how do the dynamics of model \eqref{P} change
as the advection rates vary?

Motivated by their problems, in this paper,
we shall settle their problems partially,
where we assume that both $r(x)$ and $K(x)$ are not constants.
Our main purpose is to investigate the effects of diffusion rates $d_1 ,d_2$,
the unidirectional flow $q_1,q_2$ as well as the effect of spatial heterogeneity
on the dynamical behaviors of system \eqref{P}.

For simplicity, introducing the following changes of variables
$v=aP,b=ea$,
we arrive at the following system
\beq
\label{system}
\left\{\arraycolsep=1.5pt
\begin{array}{lll}
u_t=d_1  u_{xx}-q_1 u_x+u(r(x)(1-\F{u}{K(x)})-v),&x\in(0,L),t>0,\\[2mm]
v_t=d_2  v_{xx}-q_2 v_x+v(bu-\ga),&x\in(0,L),t>0,\\[2mm]
d_1 u_x(0,t)-q_1 u(0,t)=0,\q u_x(L,t)=0,&t>0, \\[2mm]
d_2 v_x(0,t)-q_2 v(0,t)=0,\q v_x(L,t)=0,&t>0,\\[2mm]
u(x,0)=u_0(x)\ge,\not\equiv0,\q
v(x,0)=v_0(x)\ge,\not\equiv0,&x\in(0,L).
\end{array}\right.
\eeq
It seems reasonable to assume that the intrinsic growth rate $r(x)$ should have some
relationship to the carrying capacity $K(x)$;
that is, if $K(x_1)=K(x_2)$,
then $r(x_1)=r(x_2)$.
In other words, $r$ is a function of $K$,
i.e. $r(x)=h(K(x))$ for some function $h$.

From now on we assume $L=1$ for simplicity
and
\beq
\label{assumption}
\bea{l}
d_1 ,d_2 ,b,\ga,r(x),K(x)>0,q_1,q_2\ge0,r(x),\F{r(x)}{K(x)}\in C^{1}[0,1], r(x)=h(K(x)),\\
\text{both}~r(x)~\text{and}~K(x)~\text{are not constants},
\eea
\eeq
throughout this paper.

For our purpose, we will  study the relationship between
\eqref{system} and the following single species model:
\beq
\label{single}
\left\{\arraycolsep=1.5pt
\begin{array}{lll}
u_t=d u_{xx}-q u_x+r(x)u(1-\F{u}{K(x)}),&x\in(0,1),t>0,\\[2mm]
d u_x(0,t)-q u(0,t)=u_x(1,t)=0,&t>0, \\[2mm]
u(x,0)\ge,\not\equiv0,&x\in[0,1].
\end{array}\right.
\eeq
Therefore, in this paper, before studying the dynamics of system \eqref{system},
we will first reveal the dynamics of \eqref{single}.
Recall that the existence of the positive solution for \eqref{single}
depends on the following eigenvalue problem
\beq
\label{eigenvalue problem}
\left\{\arraycolsep=1.5pt
\begin{array}{lll}
d\phi_{xx}-q\phi_x+r(x)\phi=\sig\phi,&x\in(0,1),\\[2mm]
d\phi_x(0)-q\phi(0)=\phi_x(1)=0.& \\[2mm]
\end{array}\right.
\eeq
It is well known that (see, e.g., \cite{CH,HL}) problem \eqref{eigenvalue problem} admits a principal
eigenvalue $\sig_1(d,q,r)$ depending continuously and differentially on the parameters $d ,q$,
which is simple, and its corresponding eigenfunction, denoted by $\phi_1$,
can be chosen positive on [0,1].
It follows from \cite[Theorem 4.1]{LL} that for all $d>0$,
if $\sig_1(d,0,r)\le0$,
then $\sig_1(d,q,r)<0$ for any $q>0$;
if $\sig_1(d,0,r)>0$, then there exists a unique $q^*>0$
such that $\sig_1(d,q,r)>0$ for $0\le q<q^*$ and $\sig_1(d,q,r)<0$ for $q>q^*$.

Now the dynamical results on system \eqref{single} are as follows.
\begin{lem}[Section 3.2.3\cite{CC}]
\label{existence of theta}
Suppose $d>0$ fixed.
If $\sig_1(d,0,r)\le0$, then $u=0$ is globally asymptotically stable
among all solutions of \eqref{single}.
If $\sig_1(d,0,r)>0$, then there exists a unique $q^*(d,r)>0$ such that
for $q\ge q^*(d,r)$, $u=0$ is globally asymptotically stable;
and for $0<q<q^*(d,r)$,
there is a unique positive steady state of \eqref{single}
(denoted by $\th(d,q)$),
which is globally asymptotically stable.
More precisely, $\th(d,q)$ is the solution of
\beq
\label{theta}
\left\{\arraycolsep=1.5pt
\begin{array}{lll}
d\th_{xx}-q \th_x+r(x)\th(1-\F{\th}{K(x)})=0,&x\in(0,1),\\[2mm]
d\th_x(0)-q \th(0)=\th_x(1)=0.& \\[2mm]
\end{array}\right.
\eeq
\end{lem}
Clearly, $r(x)>0$ yields $\sig_1(d,0,r)>0$,
which implies that the semitrivial steady state solution $(\th_1,0)$ exists provided that $0\le q_1<q^{*}(d_1,r)$.
Here we define $\th_1=\th(d_1,q_1)$.

Notice from the similar arguments in \cite{NWW} that
each of the steady state solutions(i.e. $(0,0)$ and $(\th_1,0)$) is globally asymptotically stable if it is locally asymptotically stable, respectively.
And system \eqref{system} is uniformly persistence
if all of the steady state solutions are unstable.
Therefore, in the paper, we shall focus on the local stability of the steady state solutions.

To describe our main result,
set
\beq
\label{kappa}
\kappa_1:=\I_{0}^{1}K(x)\dx,
~\kappa_2:=\F{\i_{0}^{1}r(x)\dx}{\i_{0}^{1}\f{r(x)}{K(x)}\dx},
~\kappa_3:=\sup\limits_{d_1>0}\I_{0}^{1}\eta(x)\dx,
~\kappa_4:=\max\limits_{x\in[0,1]}K(x),
\eeq
where
$\eta$ is the unique positive solution of
\beq
\label{eta}
\left\{\arraycolsep=1.5pt
\begin{array}{lll}
d_1 \eta_{xx}+r(x)\eta(1-\F{\eta}{K(x)})=0,&x\in(0,1),\\[2mm]
\eta_{x}(0)=\eta_x(1)=0.&
\end{array}\right.
\eeq
Obviously, the assumption \eqref{assumption} leads to $\sig_1(d_1,0,r)>0$.
It follows from Lemma \ref{properties of eta} that
$$\kappa_1<\kappa_2<\kappa_3<\kappa_4$$
under some conditions.
The main results of this paper are as follows.

\begin{thm}
\label{thm1}
Assume that \eqref{assumption} holds,
$r'\ge(\F{r}{K})'\max\limits_{x\in[0,1]}K(x)>0$,
$\F{h(s)}{s}$ is strictly increasing in $s>0$
and $q_1=q_2=q$.
Let $q^{*}(d_1,r)$ be defined in Lemma \ref{existence of theta}.
Then there exists a unique $q_0(d_1,d_2,\ga)\in[0,q^{*}(d_1,r))$ such that the following conclusions
hold:

(i) If $0<\F{\ga}{b}\le\kappa_1$,
then for every $0\le q<q_0(d_1,d_2,\ga)$,
$(\th_1,0)$ is unstable for all $d_1,d_2>0$;
for every $q_0(d_1,d_2,\ga)<q<q^{*}(d_1,r)$,
$(\th_1,0)$ is stable for any $d_1,d_2>0$.

(ii) If $\kappa_1<\F{\ga}{b}<\kappa_2$, there exists a unique $d_{2*}>0$ such that

\ \ \ (a) if $0<d_2<d_{2*}$, then for
every $0\le q<q_0(d_1,d_2,\ga)$, $(\th_1,0)$ is unstable for any $d_1>0$; for every $q_0(d_1,d_2,\ga)<q<q^*(d_1,r)$, $(\th_1,0)$ is stable for any $d_1>0$;

\ \ \ (b) if $d_2>d_{2*}$,
there exists a unique $\underline{d}_1>0$
and a $\eps_0>0$
such that $(\th_1,0)$ is unstable for any $d_1>\underline{d}_1$ and $0\le q<q_0(d_1,d_2,\ga)$;
$(\th_1,0)$ is stable for any $d_1>\underline{d}_1$ and $q_0(d_1,d_2,\ga)<q<q^{*}(d_1,r)$
or $\underline{d}_1-\eps_0<d_1<\underline{d}_1$ and $q\ge0$.

(iii) If $\kappa_2<\F{\ga}{b}<\kappa_3$,
then there exists a unique $\bar{d_2}>0$ such that

\ \ \ (a) there exists a unique $\hat{d_1}>0$ such that
$(\th_1,0)$ is stable for any $d_1\ge\hat{d_1}$, $0\le q<q^*(d_1,r)$, $d_2>0$;

\ \ \ (b) there exists a $\bar{d_1}\in[0,\hat{d_1})$, such that $(\th_1,0)$ is unstable for any $d_1\in(\bar{d_1},\hat{d_1})$, $d_2<\bar{d_2}$, $0\le q<q_0(d_1,d_2,\ga)$; $(\th_1,0)$ is stable for any $d_1\in(\bar{d_1},\hat{d_1})$, $d_2<\bar{d_2}$, $q_0(d_1,d_2,\ga)<q<q^*(d_1,r)$;

\ \ \ (c) there exists a $\check{d_1}\in(0,\hat{d_1})$, such that $(\th_1,0)$ is unstable for any $d_1\in(\check{d_1},\hat{d_1})$, $d_2>\bar{d_2}$, $0\le q<q_0(d_1,d_2,\ga)$; $(\th_1,0)$ is stable for any $d_1\in(\check{d_1},\hat{d_1})$, $d_2>\bar{d_2}$, $q_0(d_1,d_2,\ga)<q<q^*(d_1,r)$.

(iv) If $\kappa_3<\F{\ga}{b}<\kappa_4$,
then there exist unique $\hat{d_1}>0$, $\tilde{d_2}>0$, $\lam^{**}>0$ such that

\ \ \ (a) $(\th_1,0)$ is stable for any $d_1\ge\hat{d_1}$,
$0\le q<q^*(d_1,r)$, $d_2>0$;
or $d_1>0$, $0\le q<q^*(d_1,r)$, $d_2>\tilde{d_2}$;

\ \ \ (b) there exists a $\bar{d_1}\in[0,\hat{d_1})$, such that $(\th_1,0)$ is unstable for any $d_1\in(\bar{d_1},\hat{d_1})$, $d_2\le\F{1}{\lam^{**}}$, $0\le q<q_0(d_1,d_2,\ga)$; $(\th_1,0)$ is stable for any $d_1\in(\bar{d_1},\hat{d_1})$, $d_2\le\F{1}{\lam^{**}}$, $q_0(d_1,d_2,\ga)<q<q^*(d_1,r)$;

\ \ \ (c) there exists a $\check{d_1}\in(0,\hat{d_1})$, such that $(\th_1,0)$ is unstable for any $d_1\in(\check{d_1},\hat{d_1})$, $\F{1}{\lam^{**}}<d_2<\tilde{d_2}$, $0\le q<q_0(d_1,d_2,\ga)$; $(\th_1,0)$ is stable for any $d_1\in(\check{d_1},\hat{d_1})$, $\F{1}{\lam^{**}}<d_2<\tilde{d_2}$, $q_0(d_1,d_2,\ga)<q<q^*(d_1,r)$.

(v) If $\ga\ge\ga_4$,
$(\th_1,0)$ is stable for any $d_1,d_2>0$, $0\le q<q^*(d_1,r)$.

(vi)
\beq
\begin{cases}
\lim\limits_{(d_1,d_2)\rightarrow(+\infty,+\infty)}q_0(d_1,d_2,\ga)
=\F{b\i_0^1r(x)\dx-\ga\i_0^1\frac{r(x)}{K(x)}\dx}{b+\i_0^1\frac{r(x)}{K(x)}\dx},\\
\lim\limits_{d_1\rightarrow+\infty}q_0(d_1,d_2,\ga)=\hat{q}_0,\\
\lim\limits_{d_2\rightarrow+\infty}q_0(d_1,d_2,\ga)=\tilde{q}_0,\\
\F{\pa}{\pa\ga}q_0(d_1,d_2,\ga)<0,
\end{cases}
\eeq
where
$\hat{q}_0$ and $\tilde{q}_0$ satisfy
\beq
\begin{cases}
\sig_1(d_2,\hat{q}_0,\F{b(\i_{0}^{1}r(x)\dx-\hat{q}_0)}{\i_{0}^{1}\f{r(x)}{K(x)}\dx}-\ga)=0,\\
b\I_0^1\th_1|_{q=\tilde{q}_0}\dx=\ga+\tilde{q}_0.
\end{cases}
\eeq
\end{thm}

We remark that Theorems \ref{thm1} indicates
the effects of diffusion $d_1 ,d_2 $, mortality rate $\ga$, advection rate $q$
on the dynamics of system \eqref{system}.
In fact, the dynamics of the predator-prey system in open advective environments
can dramatically change as the mortality rate and the advection rate vary.
\begin{figure}[htbp]
\centering
\includegraphics[width=7cm]{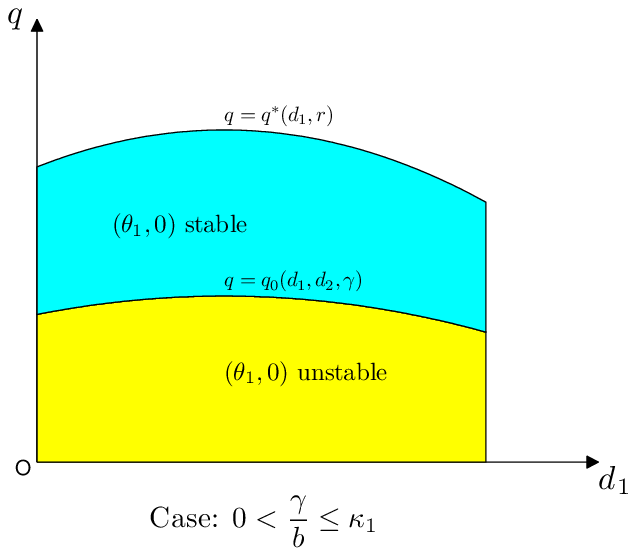}
\caption{{\scriptsize This picture illustrates part(i) of Theorems \ref{thm1},
where $\F{\ga}{b}\in(0,\kappa_1]$.
The yellow region indicates that $(\th_1,0)$ is unstable,
which means that the predator can invade when rare.
In blue region, $(\th_1,0)$ is stable
and the predator can not invade when rare.}}
\label{fig:1}
\end{figure}
\begin{figure}[htbp]
\centering
\includegraphics[width=7cm]{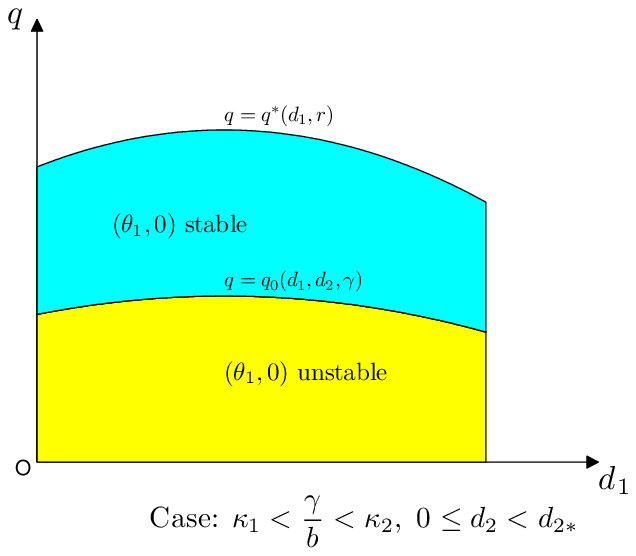}
\caption{{\scriptsize
The stability of $(\th_1,0)$ in the case
where $\F{\ga}{b}\in(\kappa_1,\kappa_2),d_2<d_{2*}$.
In the yellow region, $(\th_1,0)$ is unstable;
in the blue region, $(\th_1,0)$ is stable.
}}
\label{fig:21}
\end{figure}
\begin{figure}[htbp]
\centering
\includegraphics[width=7cm]{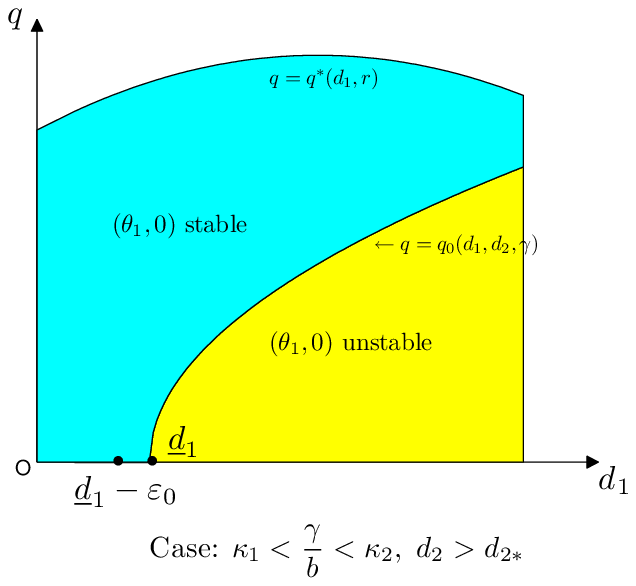}
\caption{{\scriptsize  The stability of $(\th_1,0)$ in the case
where $\F{\ga}{b}\in(\kappa_1,\kappa_2),d_2>d_{2*}$.
In the yellow region, $(\th_1,0)$ is unstable;
in the blue region, $(\th_1,0)$ is stable.
In this case, one can see from \eqref{f} that
$\sig_1(d_2,0,b\eta-\ga)$ changes sign at least once
as $d_1$ varies from $0$ to $+\infty$.
For illustration, we just assume in this figure that $\sig_1(d_2,0,b\eta-\ga)$ changes sign
exactly once in $(0,+\infty)$.}}
\label{fig:22}
\end{figure}
\begin{figure}[htbp]
\centering
\includegraphics[width=7cm]{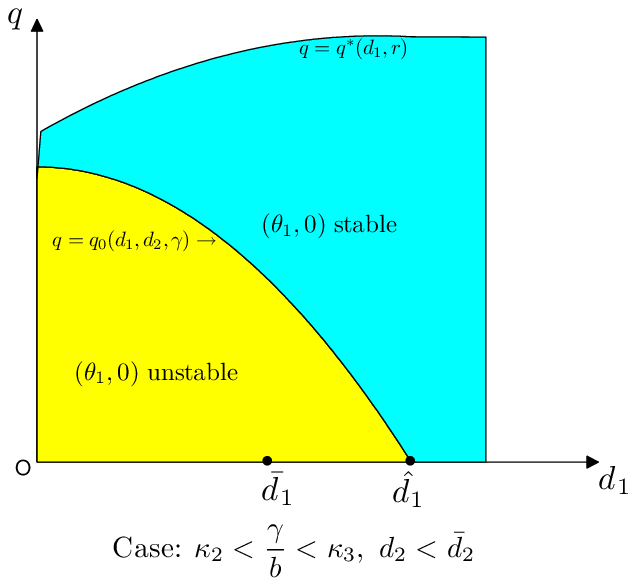}
\caption{{\scriptsize  The stability of $(\th_1,0)$ in the case
where $\F{\ga}{b}\in(\kappa_2,\kappa_3),d_2<\bar{d}_2$.
In the yellow region, $(\th_1,0)$ is unstable;
in the blue region, $(\th_1,0)$ is stable.
In this case, one can get from \eqref{f} that
$\sig_1(d_2,0,b\eta-\ga)$ changes sign at least once
as $d_1$ varies from $0$ to $+\infty$.
For illustration, we just assume in this figure that $\sig_1(d_2,0,b\eta-\ga)$ changes sign
exactly once in $(0,+\infty)$.}}
\label{fig:31}
\end{figure}
\begin{figure}[htbp]
\centering
\includegraphics[width=7cm]{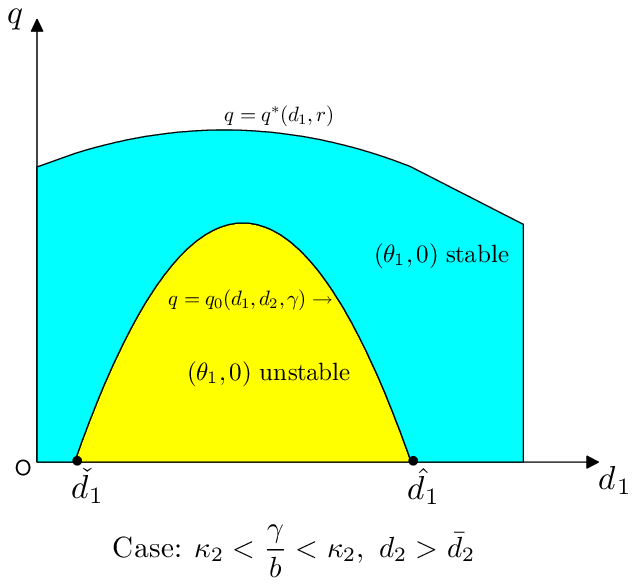}
\caption{{\scriptsize  The stability of $(\th_1,0)$ in the case
where $\F{\ga}{b}\in(\kappa_2,\kappa_3),d_2>\bar{d}_{2}$.
In the yellow region, $(\th_1,0)$ is unstable;
in the blue region, $(\th_1,0)$ is stable.
In this case, one can see that
$\sig_1(d_2,0,b\eta-\ga)$ changes sign at least twice
if $d_1\in(0,+\infty)$.
For illustration, we just assume in this figure that $\sig_1(d_2,0,b\eta-\ga)$ changes sign
exactly twice in $(0,+\infty)$.}}
\label{fig:32}
\end{figure}
\begin{figure}[htbp]
\centering
\includegraphics[width=7cm]{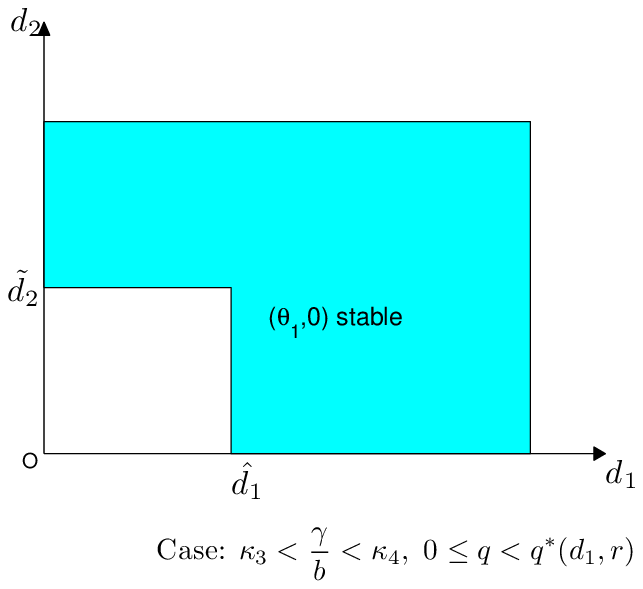}
\caption{{\scriptsize  The stability of $(\th_1,0)$ in the case
where $\F{\ga}{b}\in(\kappa_3,\kappa_4),0\le q<q^{*}(d_1,r)$.
In the blue region, $(\th_1,0)$ is stable.}}
\label{fig:41}
\end{figure}
\begin{figure}[htbp]
\centering
\includegraphics[width=7cm]{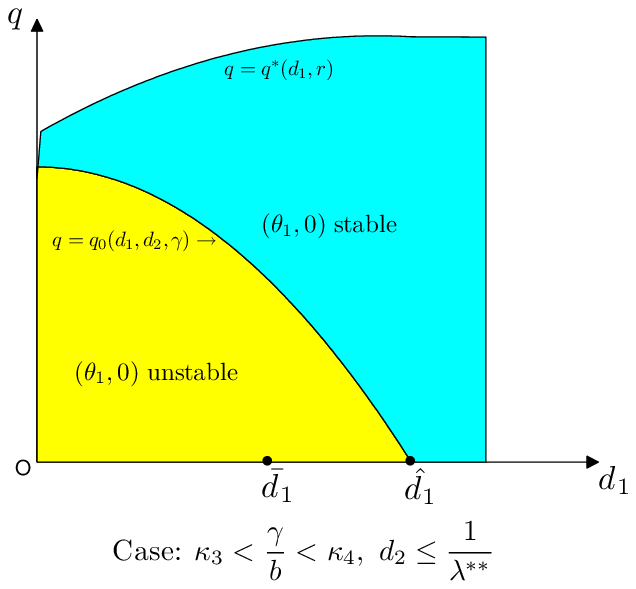}
\caption{{\scriptsize  The stability of $(\th_1,0)$ in the case
where $\F{\ga}{b}\in(\kappa_3,\kappa_4),d_2\le\F{1}{\lam^{**}}$.
In the yellow region, $(\th_1,0)$ is unstable;
in the blue region, $(\th_1,0)$ is stable.
Also in this case, one can see that
$\sig_1(d_2,0,b\eta-\ga)$ changes sign at least once
if $d_1\in(0,+\infty)$.
For illustration,
we just assume in this figure that
$\sig_1(d_2,0,b\eta-\ga)$ changes sign
exactly once for $d_1\in(0,+\infty)$.}}
\label{fig:42}
\end{figure}
\begin{figure}[htbp]
\centering
\includegraphics[width=7cm]{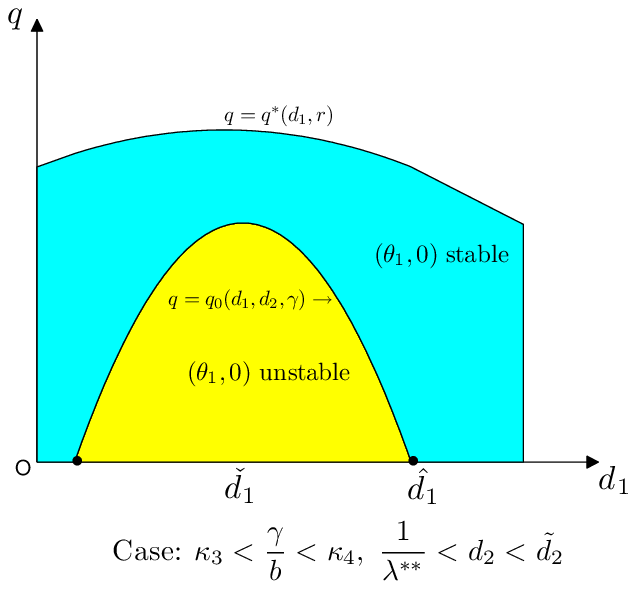}
\caption{{\scriptsize
The stability of $(\th_1,0)$ in the case
where $\F{\ga}{b}\in(\kappa_3,\kappa_4), \F{1}{\lam^{**}}<d_2<\hat{d_2}$.
In the yellow region, $(\th_1,0)$ is unstable;
in the blue region, $(\th_1,0)$ is stable.
Also in this case, one can see that
$\sig_1(d_2,0,b\eta-\ga)$ changes sign at least twice
if $d_1\in(0,+\infty)$.
For illustration, we just assume in this figure that $\sig_1(d_2,0,b\eta-\ga)$ changes sign
exactly twice for $d_1\in(0,+\infty)$.}}
\label{fig:43}
\end{figure}

Some of the conclusions of part (i),(ii), (iii),
and (iv) are illustrated in
Figs. \ref{fig:1}, \ref{fig:21}, \ref{fig:22},
\ref{fig:31}, \ref{fig:32},
\ref{fig:41}, \ref{fig:42} and \ref{fig:43}, respectively.
Note that due to the complexity of the curves separating the invasion and no-invasion regions,
Figs. \ref{fig:1}, \ref{fig:21}, \ref{fig:22},
\ref{fig:31}, \ref{fig:32}, \ref{fig:41}, \ref{fig:42}
and \ref{fig:43} are for illustration purpose only.

Biologically, the first part of part (i) means that if the death rate of the predator is
small,
then there exists a unique curve $q=q_0(d_1,d_2,\ga)$ in the $d_1-q$ plane
such that if the advection rate $q$ is smaller than $q_0(d_1,d_2,\ga)$,
then the predator can always invade when rare. In contrast, the second part of part (i) and part (v) imply that if the death rate of
the predator is large enough or
if the advection rate $q\in(q_0(d_1,d_2,\ga),q^{*}(d_1,r))$ and the death rate of the predator is
small, it can not invade.

Part (ii) is illustrated in
Fig. \ref{fig:21}, \ref{fig:22}.
More precisely, when $\kappa_1<\F{\ga}{b}<\kappa_2$,
the predator can invade when rare,
provided that one of the following two conditions holds:
\ding{172} $0\le q<q_0(d_1,d_2,\ga)$ and $d_2<d_{2*}$;
\ding{173} $0\le q<q_0(d_1,d_2,\ga)$, $d_1>\underline{d}_1$
and $d_2>d_{2*}$.

Part (iii) is illustrated in Fig. \ref{fig:31}, \ref{fig:32}.
In this part, it can be advantageous for prey to have large diffusion rate:
If the diffusion rate of the prey $d_1$ is larger than $\hat{d_1}$,
where $\hat{d_1}$ is the largest
positive root of $\max\limits_{x\in[0,1]}\eta=\ga$,
the predator can never invade when rare;
if the diffusion rate of the prey is smaller than $\hat{d_1}$,
the predator can invade when rare
by adopting smaller advection rate $q$.
Therefore, only if the diffusion rate of the prey is smaller than $\hat{d_1}$,
it may be advantageous for the predator to adopt
smaller advection rate to invade.

Part (iv) is illustrated in Fig. \ref{fig:41}, \ref{fig:42}, \ref{fig:43}.
Large diffusion rate($d_1>\hat{d_1}$) is good for the prey as the predator can never invade when rare. In contrast, if the diffusion rate of the prey is less than $\hat{d_1}$, the predator can
always adopt small dispersal rate to invade. Notice that for part (iv), if the dispersal rate of
the predator is greater than $\tilde{d_2}$ , it can not invade for any dispersal rate of the
prey.

\begin{thm}
\label{thm2}
Assume that \eqref{assumption} holds,
$r'\ge(\F{r}{K})'\max\limits_{x\in[0,1]}K(x)>0$.
Let $q^{*}(d_1,r)$ be defined in Lemma \ref{existence of theta}
and $\mu_1(q_1,q_2)$ be defined in \eqref{mu1}.
Then the following conclusions hold:

(i) $(0,0)$ is globally asymptotically stable
provided that $q_1\ge q^{*}(d_1,r)$;

(ii) if $\mu_1(0,0)\le0$,
$(q_1,q_2)\in[0,q^{*}(d_1,r))\times[0,+\infty)$,
then $(\th_1,0)$ is stable;

(iii) if $\mu_1(0,0)>0$,
$(q_1,q_2)\in\big\{(q_1,q_2)|0\le q_1<q^{*}(d_1,r),q_2\ge\tilde{q}_2\big\}
\bigcup\big\{(q_1,q_2)|\tilde{q}_1(q_2)< q_1<q^{*}(d_1,r), 0\le q_2<\tilde{q}_2\big\}$,
then $(\th_1,0)$ is stable;

(iv)if $\mu_1(0,0)>0$,
$(q_1,q_2)\in\{(q_1,q_2)|0\le q_1<\tilde{q}_1(q_2),0\le q_2<\tilde{q}_2\}$,
then $(\th_1,0)$ is unstable.

Here $\tilde{q}_2$ and $\tilde{q}_1(q_2)$ are defined in Section \ref{Proof of Theorem thm2},
and satisfy $\tilde{q}_1(\tilde{q}_2)=0$,
$\tilde{q}_1(0)=q^{*}(d_1,r)$, $\F{\mathrm{d}\tilde{q}_1(q_2)}{\mathrm{d} q_2}<0$.
\end{thm}

\begin{figure}[htbp]
\centering
\includegraphics[width=7cm]{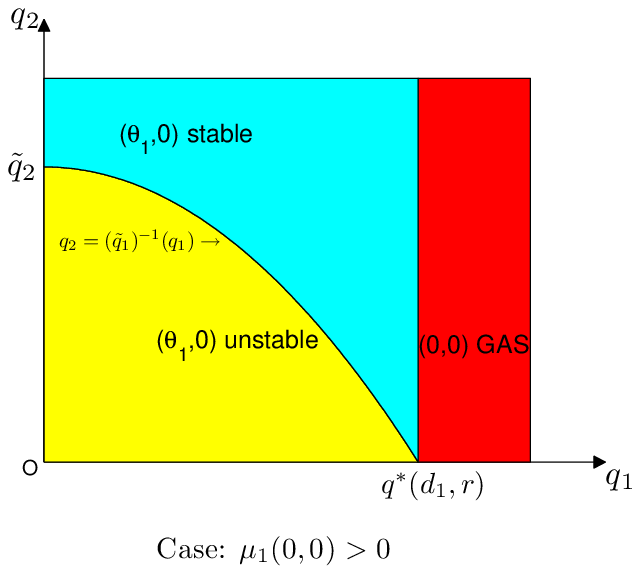}
\caption{{\scriptsize
In the yellow region, $(\th_1,0)$ is unstable;
in the blue region, $(\th_1,0)$ is stable;
in the red region, $(0,0)$ is globally asymptotically stable.}}
\label{fig:II}
\end{figure}

Theorems \ref{thm2} indicates that the dynamical behaviors of system \eqref{system}
depend heavily on the sign of $\mu_1(0,0)$,
which is similar to that of Theorem \ref{thm1}.
More precisely,
if $\mu_1(0,0)\le0$, $(\th_1,0)$ is always stable.
If $\mu_1(0,0)>0$,
there is a critical advection rate $\tilde{q}_2$ of the predator
and a critical curves $\tilde{q}_1(q_2)$ in the $q_1-q_2$ plane,
which classify the dynamic behaviors of system \eqref{system} into three scenarios
(see Fig. \ref{fig:II}).

\begin{thm}
\label{thm3}
Let $q^{*}(d_1,r)$ be defined in Lemma \ref{existence of theta}.
Assume that \eqref{assumption} holds,
$r'\ge(\F{r}{K})'\max\limits_{x\in[0,1]}K(x)>0$,
$q_1=\tau q_2$.
Then

(i) $(0,0)$ is globally asymptotically stable provided that $q_2\ge\F{q^{*}(d_1,r)}{\tau}$;

(ii) if $0\le q_2<\F{q^{*}(d_1,r)}{\tau}$,
there exists a critical curve $b=b_{\tau}(q_2)\in(0,+\infty)$
continuously depending on the parameter $q_2$
such that $(\th_1,0)$ is locally
asymptotically stable if $b\in(0,b_{\tau}(q_2))$,
and unstable if $b\in(b_{\tau}(q_2),+\infty)$.
Moreover,
the critical curve $b=b_{\tau}(q_2)\in(0,+\infty)$ is strictly increasing with respect to $q_2$,
and
\beq
\begin{cases}
b_{\tau}(0)=\inf\limits_{\om\in H^{1}([0,1])\backslash\{0\}}
\F{d_2\i_{0}^{1}\om^2_x\dx+\ga\i_0^1\om^2\dx}{\i_{0}^{1}\eta\om^2\dx}\ge\F{\ga}{\max\limits_{x\in[0,1]}K},\\
\lim\limits_{q_2\rightarrow(\frac{q^{*}(d_1,r)}{\tau})^{-}}b_{\tau}(q_2)=+\infty.
\end{cases}
\eeq
\end{thm}

\begin{figure}[htbp]
\centering
\includegraphics[width=7cm]{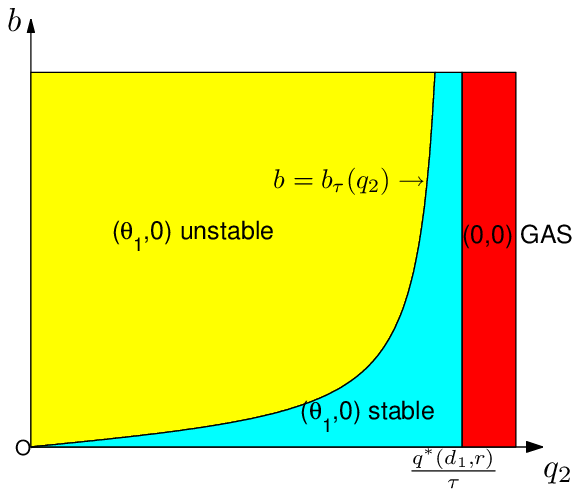}
\caption{{\scriptsize
Schematic illustrations of Theorems \ref{thm3}
with $q_1=\tau q_2$ in the $q_2-b$ plane.
Each colored region indicates is similar to Fig. \ref{fig:II}.}}
\label{fig:III}
\end{figure}

Theorems \ref{thm3} indicates that the critical curves $b=b_{\tau}(q_2)$ divide
the dynamics of system \eqref{system} with $q_1=\tau q_2$ in the $q_2-b$ plane into three scenarios (see Figs. \ref{fig:III}).
Biological speaking,
both the prey and predators will be washed out
if the predator's flow speed is large enough ($q_2\ge\F{q^{*}(d_1,r)}{\tau}$).
If the predator's flow speed is small ($q_2<\F{q^{*}(d_1,r)}{\tau}$),
then predators can invade with suitable the small predation rate,
followed by predators going extinct
when the predation rate continues to increase ($b>b_{\tau}(q_2)$).
That is, large predation rate is beneficial to predators,
if the predator's flow speed is small.

The rest of this paper is organized as follows.
In Sect. \ref{preliminaries},
we present some lemmas which shall be utilized in the subsequent analysis.
In Sect. \ref{steady state solutions},
we verify Theorem \ref{thm1}-\ref{thm3}.

\section{Preliminaries}
\label{preliminaries}

In this section we present some existing results concerning the steady state of single-species model as well as some useful lemmas in regard to the related linear eigenvalue problem as preliminaries.
With this in mind, we first recall that the following
result,
which characterizes the properties of $\eta$.

\begin{lem}[\cite{DNZ,GHN,LiLou}]
\label{properties of eta}
Assume that \eqref{assumption} holds.
Then the
following statements hold for the positive solution $\eta$ of \eqref{eta}:

(i) $\eta\rightarrow K$ in $L^{\infty}(0,1)$
as $d_1\rightarrow0+$.

(ii) $\eta\rightarrow\kappa_2$ in $L^{\infty}(0,1)$ as $d_1\rightarrow+\infty$.

(iii) If $\F{h(s)}{s}$ is (strictly) decreasing in $s>0$,
then $\kappa_2(<)\le\kappa_1$.

(iv) If $\F{h(s)}{s}$ is (strictly) increasing in $s>0$,
then $\kappa_2(>)\ge\kappa_1$.

(v) If $K(x)\not\equiv$ constant, then $\inf\limits_{x\in[0,1]}K(x)<\eta(x)<\sup\limits_{x\in[0,1]}K(x)$.


Here $\kappa_1$ and $\kappa_2$ are defined in \eqref{kappa}.
\end{lem}

\begin{lem}
\label{int eta>kappa1}
Assume that \eqref{assumption} holds and $r'\ge(\F{r}{K})'\max\limits_{x\in[0,1]}K(x)>0$.
Then $\eta_x>0$ in $(0,1)$ and $\I_0^1\eta\dx>\kappa_1$.
\end{lem}

\proof
Dividing \eqref{eta} by $\F{\eta r}{K}$ and integrating on $[0,1]$,
we obtain that
\beq
\I_0^1\eta\dx=\I_0^1 K(x)\dx+d_1\I_0^1\F{K(x)(\eta_x)^2}{r(x)\eta^2}\dx-d_1\I_0^1\F{1}{\eta}\eta_x\big(\F{K}{r}\big)_x\dx.
\eeq
By a transformation $m=\F{\eta_x}{\eta}$ and some straightforward computations, we will find
that $m$ satisfies
\beq
\left\{
\bea{ll}
d_1m_{xx}+2d_1mm_x-\F{r\eta}{K}m=(\F{r}{K})'\eta-r',&0<x<1,\\
m(0)=m(1)=0.&
\eea
\right.
\eeq
Then it would follow from the maximum principle and $r'\ge(\F{r}{K})'\max\limits_{x\in[0,1]}K(x)>0$ that $\eta_x>0$ and $\eta_x\big(\F{K}{r}\big)_x<0$.
Therefore $\I_0^1\eta\dx>\kappa_1$.
\proofend

Next consider an auxiliary linear eigenvalue problem.
\beq
\label{linear eigenvalue problem}
\left\{\arraycolsep=1.5pt
\begin{array}{lll}
d\phi_{xx}-q \phi_x+h(x)\phi=\sig\phi,&x\in(0,1),\\[2mm]
d\phi_x(0)-q \phi(0)=\phi_x(1)=0,& \\[2mm]
\end{array}\right.
\eeq
where $d>0$, $q\ge0$ and $h(x)\in C[0,1]$.
The Krein-Rutman theorem (\cite{S}) guarantees that
the problem \eqref{linear eigenvalue problem} admits a principal eigenvalue,
denoted by $\sig_1(d,q,h)$,
which corresponds to a positive eigenfunction $\phi_1$ normalized by $\max\limits_{x\in[0,1]}\phi_1=1$.
Furthermore,
by the variational characterization of the principal eigenvalue,
we have
\beq
\label{sig(dqh)}
\sig_1(d,q,h)=\sup\limits_{\om\in H^{1}([0,1])\backslash\{0\}}
\F{\i_{0}^{1}e^{-\f{q}{d}x}(-d\om^2_x+h(x)\om^2)\dx-q\om^2(0)}{\i_{0}^{1}e^{-\f{q}{d}x}\om^2\dx}.
\eeq
If $h(x)\in C[0,1]$ is not a constant,
the principal eigen-pair $(\sig_1(d,q,h),\phi_1)$ has the following properties.

\begin{lem}
\label{lemma of eigenvalue}
The following statements on the principal eigen-pair $(\sig_1(d,q,h),\phi_1)$ are true:

(i) $\sig_1(d,q,h)$ and $\phi_1(d,q,h)$ depend continuously and differentially on
parameters $d$ and $q$;

(ii) $h_n(x)\rightarrow h(x)$ in $C[0,1]$ implies $\sig_1(d,q,h_n)\rightarrow\sig_1(d,q,h)$;

(iii) $h_1(x)\ge h_2(x)$ implies $\sig_1(d,q,h_1)\ge \sig_1(d,q,h_2)$,
and the equality holds only if $h_1(x)\equiv h_2(x)$;

(iv) $\sig_1(d,q,h)$ is strictly decreasing with respect to $q$ in $[0,+\infty)$
with $\lim\limits_{q\rightarrow+\infty}\sig_1(d,q,h)=-\infty$;

(v) $\lim\limits_{d\rightarrow0^{+}}\sig_1(d,q,h)=-\infty$ if $q>0$,
and $\lim\limits_{d\rightarrow+\infty}\sig_1(d,q,h)=\I_{0}^{1}h(x)\dx-q$ if $q\ge0$;

(vi)
$\sig_1(d,0,h)$ is strictly decreasing with respect to $d$,
and
\beq
\lim\limits_{d\rightarrow0^{+}}\sig_1(d,0,h)=\max\limits_{x\in[0,1]}h(x), \lim\limits_{d\rightarrow+\infty}\sig_1(d,0,h)=\I_{0}^{1}h(x)\dx;
\eeq

(vii) the positive eigenfunction $\phi_1$ satisfies $(\phi_1)_x>0$ if $h'(x)\ge0$.
\end{lem}

\proof
The proof of statements (i)-(iii) can be obtained directly in \cite{CC}.
The proof of (iv) can refer to
\cite[Lemma 4.8, Lemma 4.9]{LL}
or \cite[Lemma 2.2]{NWW}.
The statement (v) is a direct result of \cite[Lemma 2.2]{NWW}
and statement (vi) can be deduced from \cite{N}.

It remains to show (vii).
Set $P=\F{(\phi_1)_x}{\phi_1}$.
Then it follows from \eqref{linear eigenvalue problem} that
\beq
\left\{\arraycolsep=1.5pt
\begin{array}{lll}
dP_{xx}+(2dP-q)P_x=-h'(x),&x\in(0,1),\\[2mm]
P(0)=\F{q}{d},\q P(1)=0.& \\[2mm]
\end{array}\right.
\eeq
By the strong maximum principle \cite{PW},
we get $P>0$ and $(\phi_1)_x>0$ if $h'(x)\ge0$ in $[0,1]$.
The proof is finished.
\proofend

Lemma \ref{existence of theta} guarantees
the existence of the unique positive steady state solution $\th(d,q)$
of the system \eqref{single}.
Furthermore, it is easy to see that this unique positive
steady state $\th(d,q)$ satisfies the following properties.

\begin{lem}
\label{lemma of theta}
Suppose that system \eqref{single} possesses a unique positive steady state $\th(d,q)$,
and suppose that $r(x)>0,K(x)>0$ are not constants.
Then

(i) $0<\th(d,q)<\max\limits_{x\in[0,1]}K(x)$;

(ii) $\th(d,q)$ is continuously differentiable for $d>0$,
and $\lim\limits_{d\rightarrow+\infty}\th(d,q)=\F{\i_{0}^{1}r(x)\dx-q}{\i_{0}^{1}\f{r(x)}{K(x)}\dx}$
uniformly on $[0,1]$;

(iii) $\th(d,q)$ is continuously differentiable for $q$,
it is decreasing pointwisely on [0,1] when $q$ increases,
and $\th_x>0$ on $[0,1]$,
provided that
$r'(x)\ge(\F{r(x)}{K(x)})'\max\limits_{x\in[0,1]}K(x)>0$ on $[0,1]$.
Moreover, $\lim\limits_{q\rightarrow0^+}\th(d,q)=\eta$ uniformly on [0,1].
Here $\eta$ is the unique positive solution of \eqref{eta}.
\end{lem}

\proof
The statement (i) can refer to \cite{DNZ,GHN} and
the continuously differentiability can refer to \cite{CC}.
Integrating now the first equation of \eqref{theta} on $[0,x]$,
we have that
\beq
d\th_x-q\th+\i^{x}_{0}r\th(1-\F{\th}{K})\dx=0.
\eeq
Since $0<\th(d ,q)<\max\limits_{x\in[0,1]}K(x)$,
we get that $\th_x$ is uniformly bounded when $d\rightarrow+\infty$,
and furthermore \eqref{theta} reads that
$\th_{xx}$ is uniformly bounded when $d\rightarrow+\infty$.
Passing to a subsequence if necessary,
we can deduce that $\th$ tends to a positive constant $C$
uniformly on $[0,1]$ as $d\rightarrow+\infty$.
Integrating the first equation of \eqref{theta} on $[0,1]$,
we see that
\beq
q\th(1)=\i_{0}^{1}r(x)\th(1-\F{\th}{K})\dx,
\eeq
which implies that $\lim\limits_{d\rightarrow+\infty}
\th=\F{\i_{0}^{1}r(x)\dx-q}{\i_{0}^{1}\f{r(x)}{K(x)}\dx}$.

Finally we prove the statement (iii).
By a transformation $w=\F{\th_x}{\th}$,
and some straightforward computations,
we find that $w$ satisfies
\beq
\left\{\arraycolsep=1.5pt
\begin{array}{lll}
-d w_{xx}-(2dw-q)w_x+\xi\th w=r'-(\F{r}{K})'\th,&x\in(0,1),\\[2mm]
w(0)=\F{q}{d},\q w(1)=0.& \\[2mm]
\end{array}\right.
\eeq
It then follows from the maximum principle that $\th_x>0$ on $[0,1]$
under the condition of statement (iii).
Next differentiating the equation of $\th$ with respect to $q$,
we obtain
\beq
\left\{\arraycolsep=1.5pt
\begin{array}{lll}
-d\F{\pa\th_{xx}}{\pa q}+q \F{\pa\th_{x}}{\pa q}-\F{\pa\th}{\pa q}(r-\F{2r}{K}\th)=-\th_x,&x\in(0,1),\\[2mm]
-d\F{\pa\th_{x}}{\pa q}(0)+q \F{\pa\th}{\pa q}(0)=-\th(0)<0,\q\F{\pa\th_{x}}{\pa q}(1)=0.& \\[2mm]
\end{array}\right.
\eeq
Since $\F{r(x)}{K(x)}>0$, by lemma \ref{lemma of eigenvalue} one then can see
$\sig_1(d,q,r-\F{2r}{K}\th)<\sig_1(d,q,r-\F{r}{K}\th)=0$.
Hence the strong maximum principle gives rise to
$\F{\pa\th}{\pa q}<0$.
The proof is complete.
\proofend

In order to investigate the stability of the semitrivial steady state $(\th_1,0)$,
we need to consider the following linear eigenvalue problem
\beq
\label{linearity at (th,0)}
\left\{\arraycolsep=1.5pt
\begin{array}{lll}
d_2 \psi_{xx}-q_2\psi_x+(b\th_1-\ga)\psi=\sig\psi,&x\in(0,1),\\[2mm]
d_2 \psi_x(0)-q_2\psi(0)=\psi_x(1)=0.& \\[2mm]
\end{array}\right.
\eeq
Similarly, one can conclude that the principal eigenvalue of \eqref{linearity at (th,0)},
denoted by $\sig_1(d_2,q_2,b\th_1-\ga)$,
depends continuously and differentially on parameters $d_1,d_2,q_1,q_2,b,\ga$.

\section{Dynamics of \eqref{system}}
\label{steady state solutions}

In this section,
our goal is to establish the local stability of
$(\th_1,0)$.
Hence, we need focus on the following system.
\beq
\label{steady state system}
\left\{\arraycolsep=1.5pt
\begin{array}{lll}
d_1 u_{xx}-q_1 u_x+u(r(x)(1-\F{u}{K(x)})-v)=0,&x\in(0,1),\\[2mm]
d_2 v_{xx}-q_2 v_x+v(bu-\ga)=0,&x\in(0,1),\\[2mm]
d_1 u_x(0)-q_1u(0)=0,\q u_x(1)=0,&\\[2mm]
d_2 v_x(0)-q_2v(0)=0,\q v_x(1)=0.&
\end{array}\right.
\eeq
Here, we point out that compared to the work in \cite{NWW},
due to the spatial heterogeneity of the
intrinsic growth rate and the carrying capacity of the prey,
some different methods and arguments are needed
to obtain the stability of $(\th_1,0)$.

It is well known
that the stability of $(\th_1,0)$ is determined by
the sign of the principal eigenvalue $\sig_1(d_2,q_2,b\th_1-\ga)$ of \eqref{linearity at (th,0)}.
That is, $(\th_1,0)$ is stable if $\sig_1(d_2,q_2,b\th-\ga)<0$,
and it is unstable if $\sig_1(d_2,q_2,b\th_1-\ga)>0$.

\subsection{Case $q_1=q_2=q$}
\label{case q1=q2}

In this subsection, we shall consider a special case where $q_1=q_2=q$.
Combined with Lemma \ref{lemma of eigenvalue}(iii-iv)
and Lemma \ref{lemma of theta}(iii),
one can see that
\beq
\label{decreasing of sig}
\F{\pa}{\pa q}\sig_1(d_2,q,b\th_1-\ga)<0.
\eeq
Since by Lemma \ref{existence of theta} and Lemma \ref{lemma of eigenvalue}(ii) one can get that
\beq
\lim\limits_{q\rightarrow q^{*}(d_1,r)}
\sig_1(d_2,q,b\th_1-\ga)=\sig_1(d_2,q^{*}(d_1,r),-\ga)
<-\ga<0,
\eeq
it then follows from \eqref{decreasing of sig} that the sign of $\sig_1(d_2,q,b\th_1-\ga)$ is relevant to the sign of $\sig_1(d_2,0,b\eta-\ga)$.
More precisely, when $\sig_1(d_2,0,b\eta-\ga)\le0$, then $\sig_1(d_2,q,b\th_1-\ga)<0$ for all $0<q<q^{*}(d_1,r)$;
if $\sig_1(d_2,0,b\eta-\ga)>0$,
then there exists a unique $\tilde{q}\in(0,q^{*}(d_1,r))$ such that
$\sig_1(d_2,q,b\th_1-\ga)>0$ for $0<q<\tilde{q}$
and $\sig_1(d_2,q,b\th_1-\ga)<0$ for $\tilde{q}<q<q^{*}(d_1,r)$. Therefore, to study the local stability of $(\th_1,0)$,
it is very important to consider the sign of $\sig_1(d_2,0,b\eta-\ga)$.

\subsubsection{Sign of $\sig_1(d_2,0,b\eta-\ga)$}
\label{sign of mu(0,0)}

Clearly, from \eqref{sig(dqh)},
we have that
\beq
\sig_1(d_2,0,b\eta-\ga)=
\sup\limits_{\om\in H^{1}([0,1])\backslash\{0\}}
\F{\i_{0}^{1}(-d_2\om^2_x+(b\eta-\ga)\om^2)\dx}{\i_{0}^{1}\om^2\dx}.
\eeq

Before investigating the sign of $\sig_1(d_2,0,b\eta-\ga)$,
we shall display some interesting lemmas.

\begin{lem}[\cite{CC,N}]
\label{lam*}
If $b\eta-\ga$ is positive somewhere in $(0,1)$, then
\beq
\lam^{*}(d_1):=\inf\limits_{\{\varphi\in H^1(0,1),\i_0^1(b\eta-\ga)\varphi^2\dx>0\}}\F{\i_0^1|\n\varphi|^2\dx}{\i_0^1(b\eta-\ga)\varphi^2\dx}
\eeq
is well defined, i.e., $\lam^{*}(d_1)$ is non-negative and finite, and there exists $\psi\in C^2([0,1])$ such that
\beq
\left\{\arraycolsep=1.5pt
\begin{array}{lll}
\psi_{xx}+\lam^{*}(b\eta-\ga)\psi=0,&x\in(0,1),\\[2mm]
\psi_x(0)=\psi_x(1)=0.& \\[2mm]
\end{array}\right.
\eeq
Moreover, $\lam^{*}=0$ if $b\I_0^1\eta\dx\ge\ga$ and $\lam^{*}>0$ if $b\I_0^1\eta\dx<\ga$.
\end{lem}

\begin{lem}[\cite{CC,N}]
\label{connection sig sign}
Suppose that $b\eta-\ga$ is positive somewhere in $(0,1)$.
Then $\sig_1(d_2,0,b\eta-\ga)<0$, if $d_2>\F{1}{\lam^{*}(d_1)}$,
and $\sig_1(d_2,0,b\eta-\ga)>0$, if $d_2<\F{1}{\lam^{*}(d_1)}$.
\end{lem}

Now we shall apply and modify some arguments in \cite{LW} to study the sign of $\sig_1(d_2,0,b\eta-\ga)$.
More precisely, we will divide our investigation into five cases.

\textbf{Case $\F{\ga}{b}\ge \kappa_4$.}
In this case, it follows that for any $q\in(0,q^*(d_1,r))$,
\beq
\sig_1(d_2,0,b\eta-\ga)<\sig_1(d_2,0,b\max\limits_{x\in[0,1]}K(x)-\ga)
=b\max\limits_{x\in[0,1]}K(x)-\ga\le0.
\eeq

\textbf{Case $0<\F{\ga}{b}\le\kappa_1$.}
By Lemma \ref{int eta>kappa1}, we have $\I_{0}^{1}(b\eta-\ga)\dx>0$ in this case,
which implies that $\sig_{1}(d_2 ,0,b\eta-\ga)>0$ for all $d_2 >0$.

\textbf{Case $\F{\ga}{b}\in(\kappa_1,\kappa_2)$}.
This case implies that there exist positive constants
$0<d_{1*}\le d_1^{*}$ such that
\beq
\label{ineq1}
\left\{
\bea{l}
b\I_0^1\eta\dx<\ga~\text{for}~d_1\in(0,d_{1*}),\\
b\I_0^1\eta\dx>\ga~\text{for}~d_1>d_1^{*},
\eea
\right.
\eeq
where $d_{1*},d_1^{*}$ are the
smallest and largest positive roots of $b\I_0^1\eta\dx=\ga$, respectively.

We note that since $\F{\ga}{b}<\kappa_2<\kappa_4$,
then from Lemma \ref{properties of eta} (ii-iv),
it indicates that for every $d_1>0$, $\ga<b\max\limits_{x\in[0,1]}\eta<b\max\limits_{x\in[0,1]}K(x)$ and we obtain that $b\eta-\ga$, $bK-\ga$ are always positive somewhere in $(0,1)$. Therefore $\lam^{*}(d_1)$ and $\lim\limits_{d_1\rightarrow0}\lam^{*}(d_1)$ are well defined.

By \eqref{ineq1}, $\sig_{1}(d_2,0,b\eta-\ga)>0$ holds for every
$d_1\ge d_1^{*}$.
Moreover, Lemma \ref{lam*} reads that $\lam^{*}(d_1)>0$ for
$d_1\in(0,d_{1*})$ and $\lam^{*}(d_1)=0$ for
$d_1\ge d_1^{*}$. Therefore $\lam^{*}(d_1)$ is not identically zero in $[0,d_1^{*}]$, and vanishes in $[d_1^{*},+\infty)$.

Define now
\beq
d_{2*}:=\F{1}{\sup\limits_{d_1>0}\lam^{*}(d_1)}>0.
\eeq
Consequently, when $d_2<d_{2*}$,
we have $\lam^{*}(d_1)<\F{1}{d_2}$. It then follows from Lemma \ref{lam*} that $\sig(d_2,0,b\eta-\ga)>0$.

When $d_2>d_{2*}$, since $\lam^{*}(d_1)=0$ for $d_1\ge d_1^{*}$ and $\lim\limits_{d_1\rightarrow0}\lam^{*}(d_1)$ exists, then $\sup\limits_{d_1>0}\lam^{*}(d_1)=\lam^{*}(\bar{d_1})>\F{1}{d_2}$ for some $\bar{d_1}\in[0,d_1^{*}]$. Hence $\lam^{*}(d_1)-\F{1}{d_2}$ changes sign at least once in $(0,d_1^{*})$. By
Lemma \ref{connection sig sign}, $\sig_1(d_2,0,b\eta-\ga)$ also changes sign at least once in $(0,d_1^{*})$.
Furthermore, there exists a unique $\underline{d}_1\in(0,d_1^*)$ such that
$\sig_{1}(d_2,0,b\eta-\ga)>0$ for every
$d_1>\underline{d}_1$;
there exists a $\eps_0>0$ such that $\sig_{1}(d_2,0,b\eta-\ga)<0$ for $\underline{d}_1-\eps_0<d_1<\underline{d}_1$.

\textbf{Case $\F{\ga}{b}\in(\kappa_2,\kappa_3)$.}
For this range,
similar to the case above,
there exist positive constants
$0<d_{1*}<d_1^{*}$
such that
\beq
\label{aux}
\left\{
\bea{l}
b\I_0^1\eta\dx<\ga~ \text{for}~d_1\in(0,d_{1*})\bigcup(d_1^{*},+\infty),\\
b\I_0^1\eta\dx=\ga~\text{for}~ d_1=d_{1*},d_1^{*}.
\eea
\right.
\eeq

Note from Lemma \ref{properties of eta}(i)(ii) that
\beq
b\lim\limits_{d_1\rightarrow0}\max\limits_{x\in[0,1]}\eta=b\max\limits_{x\in[0,1]}K(x)>\ga,~ b\lim\limits_{d_1\rightarrow+\infty}\max\limits_{x\in[0,1]}\eta=b\F{\i_{0}^{1}r(x)\dx}{\i_{0}^{1}\f{r(x)}{K(x)}\dx}<\ga.
\eeq
Then we get that there exists a unique $\hat{d_1}\in(0,+\infty)$ such that
\beq
\label{b}
\begin{cases}
b\max\limits_{x\in[0,1]}\eta=\ga~\text{for}~d_1=\hat{d_1},\\
\ga>b\max\limits_{x\in[0,1]}\eta~\text{for}~d_1>\hat{d_1}.
\end{cases}
\eeq
This indicates that $d_1^*<\hat{d_1}$.

\emph{Claim:}
\beq
\label{claim1}
\sig_{1}(d_2,0,b\eta-\ga)|_{d_1=d_1^*}>0,
~\sig_{1}(d_2,0,b\eta-\ga)|_{d_1=\hat{d}_1}<0.
\eeq

Clearly, \eqref{aux} yields that
\beq
\sig_1(d_2,0,b\eta-\ga)|_{d_1=d_1^*}=
\sup\limits_{\om\in H^{1}([0,1])\backslash\{0\}}
\F{\i_{0}^{1}(-d_2\om^2_x+(b\eta-\ga)\om^2)\dx}{\i_{0}^{1}\om^2\dx}
>\I_0^1(b\eta-\ga)\dx=0.
\eeq
When $d_1=\hat{d_1}$,
it indicates from Lemma \ref{lemma of eigenvalue}(iii) that
\beq
\sig_1(d_2,0,b\eta-\ga)|_{d_1=\hat{d_1}}<\sig_1(d_2,0,b\max\limits_{x\in[0,1]}\eta-\ga)|_{d_1=\hat{d_1}}=\sig_1(d_2,0,0)=0.
\eeq
Hence Claim \eqref{claim1} is correct,
and furthermore the second line of \eqref{b} yields that
\beq
\label{d}
\sig_1(d_2,0,b\eta-\ga)<0~\text{for all}~ d_1\ge\hat{d_1}.
\eeq

Notice that
\beq
\label{a0}
\begin{cases}
\lim\limits_{d_1\rightarrow0^+}(b\max\limits_{x\in[0,1]}\eta-\ga)=b\max\limits_{x\in[0,1]}K-\ga>0,\\
\lim\limits_{(d_1,d_2)\rightarrow(0^+,0^+)}\sig_1(d_2,0,b\eta-\ga)=\lim\limits_{d_1\rightarrow0^+}\max\limits_{x\in[0,1]}(b\eta-\ga)=b\max\limits_{x\in[0,1]}K-\ga>0,\\
\lim\limits_{d_2\rightarrow+\infty}\sig_1(d_2,0,b\eta-\ga)=b\I_0^1\eta\dx-\ga<0~\text{for}~d_1\in(0,d_{1*}).
\end{cases}
\eeq
Then we can choose $d_{10}\in(0,d_{1*})$ such that there holds
\beq
\label{a}
\lim\limits_{d_2\rightarrow0^+}\sig_1(d_2,0,b\eta-\ga)>0~\text{for any}~d_1\in(0,d_{10}].
\eeq
Therefore, for any $d_1\in(0,d_{10}]$,
due to Lemma \ref{lemma of eigenvalue}(vi),
there exists a unique $d_2(d_1)>0$ such that $\sig_1(d_2(d_1),0,b\eta-\ga)=0$.

Define \beq
\bar{d}_2:=\inf\limits_{0<d_1\le d_{10}}
\{d_2(d_1)|\sig_1(d_2,0,b\eta-\ga)=0\}.
\eeq
We claim that $\bar{d_2}$ is well-defined.
By the variational characterization of $\sig_1(d_2,0,b\eta-\ga)$,
we deduce that
\beq
\bea{l}
\sig_1(d_2,0,b\eta-\ga)=
\sup\limits_{\om\in H^{1}([0,1])\backslash\{0\}}
\F{\i_{0}^{1}(-d_2\om^2_x+(b\eta-\ga)\om^2)\dx}{\i_{0}^{1}\om^2\dx}\\
=\sup\limits_{\om\in H^{1}([0,1])\backslash\{0\}}
\F{-d_2\i_{0}^{1}\om^2_x\dx+\i_{0}^{1}(bK-\ga)\om^2\dx+\i_{0}^{1}(b\eta-bK)\om^2\dx}{\i_{0}^{1}\om^2\dx}\\
\ge\sup\limits_{\om\in H^{1}([0,1])\backslash\{0\}}
\F{-d_2\i_{0}^{1}\om^2_x\dx+\i_{0}^{1}(bK-\ga)\om^2\dx}{\i_{0}^{1}\om^2\dx}
+\inf\limits_{\om\in H^{1}([0,1])\backslash\{0\}}
\F{\i_{0}^{1}(b\eta-bK)\om^2\dx}{\i_{0}^{1}\om^2\dx}\\
:=\sig_1(d_2,0,bK-\ga)+I.
\eea
\eeq
By Lemma \ref{properties of eta}(i),
we have that there exists $d_{100}\in(0,d_{10})$ such that
\beq
I<\F{b\kappa_4-b\kappa_3}{2},~\text{for any}~d_1\in(0,d_{100}).
\eeq
Since by Lemma \ref{lemma of eigenvalue}(vi) there holds
\beq
\lim\limits_{d_2\rightarrow0^+}\sig_1(d_2,0,bK-\ga)=b\kappa_4-\ga>b\kappa_4-b\kappa_3>0,
\eeq
we then get that there exists $\delta>0$ such that
\beq
\sig_1(d_2,0,b\eta-\ga)>b\kappa_4-b\kappa_3-\F{b\kappa_4-b\kappa_3}{2}=\F{b\kappa_4-b\kappa_3}{2}>0~\text{for all}~d_2\in(0,\delta)~\text{and}~d_1\in(0,d_{100}).
\eeq
Hence, for any $d_1\in(0,d_{100})$,
if $\sig_1(d_2(d_1),0,b\eta-\ga)=0$,
then $d_2(d_1)\ge\delta>0$.
For any $d_1\in[d_{100},d_{10}]$,
the continuity of $d_2(d_1)$ with respect to $d_1$ yields that
there exists some $C\ge0$ independent of $d_1$ such that $d_2(d_1)\ge C$.
If $C=0$, one can see that $d_2(D_1)=0$ for some $D_1\in[d_{100},d_{10}]$.
However, by Lemma \ref{lemma of eigenvalue}(vi),
it indicates that $b\max\limits_{x\in[0,1]}\eta|_{d_1=D_1}-\ga=0$,
which contradicts to \eqref{a}.
Hence $C>0$. Furthermore, we obtain that $\bar{d_2}>0$.

On the other hand, for $d_1\in(0,d_{10}]\subset(0,d_{1*})$,
we have that
\beq
\lim\limits_{d_2\rightarrow+\infty}\sig_1(d_2,0,b\eta-\ga)=b\I_0^1\eta\dx-\ga<0,
\eeq
which implies that $\bar{d_2}$ is bounded above.
Therefore $\bar{d_2}$ is well-defined.

Applying Lemma \ref{lemma of eigenvalue}(vi)
with \eqref{a} and the third line of \eqref{a0},
we obtain that
if $d_2<\bar{d_2}$,
$\sig_1(d_2,0,b\eta-\ga)>0$ for any $d_1\in(0,d_{10}]$;
if $d_2>\bar{d_2}$,
there exists a $d_1\in(0,d_{10}]$ such that $\sig_1(d_2,0,b\eta-\ga)<0$.

Therefore,
when $d_2<\bar{d_2}$, there always holds $\sig_1(d_2,0,b\eta-\ga)>0$ for all $d_1\in(0,d_{10}]$.
Combined with Claim \eqref{claim1},
it follows that $\sig_1(d_2,0,b\eta-\ga)$ changes its sign at least once as $d_1$ varies from $0$ to $+\infty$.
When $d_2>\bar{d_2}$,
there exists a $d_1\in(0,d_{10}]\subset(0,d_{1*})$ such that $\sig_1(d_2,0,b\eta-\ga)<0$. Combined with Claim \eqref{claim1} again,
it follows that $\sig_1(d_2,0,b\eta-\ga)$ changes its sign at least twice as $d_1$ varies from $0$ to $+\infty$.

\textbf{Case $\F{\ga}{b}\in(\kappa_3,\kappa_4)$.}
Similarly, in this case, \eqref{d} also holds.
It is easy to see that
\beq
\label{c}
\lim\limits_{(d_1,d_2)\rightarrow(0^+,0^+)}\sig_1(d_2,0,b\eta-\ga)=\lim\limits_{d_1\rightarrow0^+}\max\limits_{x\in[0,1]}(b\eta-\ga)
=b\max\limits_{x\in[0,1]}K-\ga>0.
\eeq
Since \eqref{b} reads that, if $d_1>\hat{d_1}$, $\lim\limits_{d_2\rightarrow0^+}\sig_1(d_2,0,b\eta-\ga)=b\max\limits_{x\in[0,1]}\eta-\ga<0$,
using \eqref{c},
one then can have that there exists some $D_1\in(0,\hat{d_1})$ such that $\lim\limits_{d_2\rightarrow0^+}\sig_1(d_2,0,b\eta-\ga)>0$ for $d_1\in(0,D_1]$.
Applying $\lim\limits_{d_2\rightarrow+\infty}\sig_1(d_2,0,b\eta-\ga)=b\I_0^1\eta-\ga<0$,
we obtain that for any $d_1\in(0,D_1]$,
there exists a $d_2(d_1)>0$ such that $\sig_1(d_2(d_1),0,b\eta-\ga)=0$.

Define
\beq
\tilde{d_2}=\sup\limits_{0<d_1\le\hat{d_1}}\{d_2(d_1)|\sig_1(d_2(d_1),0,b\eta-\ga)=0\}>0.
\eeq
By variational characterization,
\beq
\bea{l}
\sig_1(d_2,0,b\eta-\ga)=
\sup\limits_{\om\in H^{1}([0,1])\backslash\{0\}}
\F{\i_{0}^{1}(-d_2\om^2_x+(b\eta-\ga)\om^2)\dx}{\i_{0}^{1}\om^2\dx}\\
=\sup\limits_{\om\in H^{1}([0,1])\backslash\{0\}}
\F{-d_2\i_{0}^{1}\om^2_x\dx+\i_{0}^{1}(bK-\ga)\om^2\dx+\i_{0}^{1}(b\eta-bK)\om^2\dx}{\i_{0}^{1}\om^2\dx}\\
\le\sup\limits_{\om\in H^{1}([0,1])\backslash\{0\}}
\F{-d_2\i_{0}^{1}\om^2_x\dx+\i_{0}^{1}(bK-\ga)\om^2\dx}{\i_{0}^{1}\om^2\dx}
+\sup\limits_{\om\in H^{1}([0,1])\backslash\{0\}}
\F{\i_{0}^{1}(bK-b\eta)\om^2\dx}{\i_{0}^{1}\om^2\dx}\\
:=\sig_1(d_2,0,bK-\ga)+II.
\eea
\eeq
Clearly, $\lim\limits_{d_1\rightarrow0^+}II=0$.
Then there exists $\underline{d_1}>0$ such that
$\sup\limits_{0<d_1<\underline{d_1}}II<\F{b(\kappa_2-\kappa_1)}{2}$.
Note that $\sig_1(d_2,0,bK-\ga)=\lim\limits_{d_1\rightarrow0^+}\sig_1(d_2,0,b\eta-\ga)$
and
$\lim\limits_{d_2\rightarrow+\infty}\sig_1(d_2,0,bK-\ga)=b\kappa_1-\ga<b(\kappa_1-\kappa_2)<0$.
Hence, there exists a $M>0$ such that for any $d_2>M$ and $d_1\in(0,\underline{d_1})$, $\sig_1(d_2,0,b\eta-\ga)<0$ for every $d_1\in(0,\hat{d_1})$,
\beq
\bea{l}
\sig_1(d_2,0,b\eta-\ga)
\le\sig_1(d_2,0,bK-\ga)+II<\F{b(\kappa_1-\kappa_2)}{2}<0.
\eea
\eeq
This inequality indicates that
if $\sig_1(d_2(d_1),0,b\eta-\ga)=0$,
then $d_2(d_1)\le M$ for every $d_1\in(0,\underline{d_1})$.
When $d_1\in[\underline{d_1},\hat{d_1}]$,
we see from the continuity of $d_2(d_1)$ with respect to $d_1$ that $d_2(d_1)\le M_1$ for some $M_1>0$.
Therefore $\tilde{d_2}$ is bounded above and furthermore is well-defined.

When $d_2<\tilde{d_2}$,
we see from \eqref{c} that there exists a $d_1\in(0,\hat{d_1}]$,
such that $\sig_1(d_2,0,b\eta-\ga)>0$.
Note from Lemma \ref{properties of eta}(i) that $\lim\limits_{d_1\rightarrow0^+}\sig_1(d_2,0,b\eta-\ga)=\sig_1(d_2,0,bK-\ga)$.
Since $\F{\ga}{b}\in(\kappa_3,\kappa_4)$,
it then follows that $bK-\ga$ is positive somewhere in $(0,1)$ and $b\I_0^1K\dx-\ga<0$.
Similar to Lemma \ref{lam*},
we have that
\beq
\lam^{**}:=\inf\limits_{\{\varphi\in H^1(0,1),\i_0^1(bK-\ga)\varphi^2\dx>0\}}\F{\i_0^1|\n\varphi|^2\dx}{\i_0^1(bK-\ga)\varphi^2\dx}>0.
\eeq
Similar to Lemma \ref{connection sig sign},
we then obtain that
\beq
\label{e}
\begin{cases}
\lim\limits_{d_1\rightarrow0^+}\sig_1(d_2,0,b\eta-\ga)=\sig_1(d_2,0,bK-\ga)<0,~\text{if}~d_2>\F{1}{\lam^{**}},\\
\lim\limits_{d_1\rightarrow0^+}\sig_1(d_2,0,b\eta-\ga)=\sig_1(d_2,0,bK-\ga)=0,~\text{if}~d_2=\F{1}{\lam^{**}},\\
\lim\limits_{d_1\rightarrow0^+}\sig_1(d_2,0,b\eta-\ga)=\sig_1(d_2,0,bK-\ga)>0,~\text{if}~d_2<\F{1}{\lam^{**}},
\end{cases}
\eeq
which implies that $\F{1}{\lam^{**}}=d_2(0)\le\tilde{d_2}$ immediately.

\emph{Claim:}
\beq
\label{claim2}
\F{1}{\lam^{**}}<\tilde{d_2}.
\eeq

Suppose on the contrary that $\F{1}{\lam^{**}}=\tilde{d_2}$.
Then we claim that there exists a $\delta>0$ such that for any $d_1\in(0,\delta)$,
there holds
$d_2(d_1)<\F{1}{\lam^{**}}$ and $\sig_1(d_2(d_1),0,b\eta-\ga)=0$.
If not, then there exists $\mathbbm{d}_1,\mathbbm{d}_2>0$ such that $\mathbbm{d}_1\neq\mathbbm{d}_2$,
$d_2(\mathbbm{d}_1)=d_2(\mathbbm{d}_2)=\F{1}{\lam^{**}}=d_2(0)$,
which implies that
\beq
\sig_1(d_2(\mathbbm{d}_2),0,b\eta|_{d_1=\mathbbm{d}_1}-\ga)=
\sig_1(d_2(\mathbbm{d}_1),0,b\eta|_{d_1=\mathbbm{d}_1}-\ga)=\sig_1(d_2(\mathbbm{d}_2),0,b\eta|_{d_1=\mathbbm{d}_2}-\ga)=0.
\eeq
Lemma \ref{lemma of eigenvalue}(iii) then yields that $\eta|_{d_1=\mathbbm{d}_1}\equiv\eta|_{d_1=\mathbbm{d}_2}:=\eta$,
which combined with \eqref{eta} implies that
\beq
\sig_1(\mathbbm{d}_1,0,r(x)(1-\F{\eta}{K(x)}))=\sig_1(\mathbbm{d}_2,0,r(x)(1-\F{\eta}{K(x)}))=0.
\eeq
Since this contradicts to Lemma \ref{lemma of eigenvalue}(vi), then $\delta$ exists.
Therefore for any $d_1\in(0,\delta)$,
there holds
$d_2(d_1)<\F{1}{\lam^{**}}$ and $\sig_1(d_2(d_1),0,b\eta-\ga)=0$.
However, by the third line of \eqref{e},
one can see that for $d_1\in(0,\delta)$ sufficiently small,
$\sig_1(d_2(d_1),0,b\eta-\ga)>0$, which is impossible.
Hence Claim \eqref{claim2} is correct.

On the other hand,
since $\F{\ga}{b}\in(\kappa_3,\kappa_4)$ leads to $b\kappa_2-\ga<0$,
it then indicates from Lemma \ref{properties of eta}(ii) and Lemma \ref{lemma of eigenvalue}(iii) that
$\lim\limits_{d_1\rightarrow+\infty}\sig_1(d_2,0,b\eta-\ga)=\sig_1(d_2,0,b\kappa_2-\ga)<0$.

Consequently, by Lemma \ref{lemma of eigenvalue}(vi),
we see that for any $d_2>\tilde{d_2}$,
$\sig_1(d_2,0,b\eta-\ga)<0$.
While $d_2\le\F{1}{\lam^{**}}$,
we obtain that $\sig_1(d_2,0,b\eta-\ga)$ changes sign at least once as $d_1$ varies from $0$ to $+\infty$.
When $\F{1}{\lam^{**}}<d_2<\tilde{d_2}$,
we obtain that $\sig_1(d_2,0,b\eta-\ga)$ changes sign at least twice as $d_1$ varies from $0$ to $+\infty$.

In summary, we can conclude that
\beq
\label{f}
\begin{cases}
\sig_{1}(d_2 ,0,b\eta-\ga)>0,~\text{if}~0<\F{\ga}{b}\le\kappa_1,\\
\sig_{1}(d_2 ,0,b\eta-\ga)>0,~\text{if}~\kappa_1<\F{\ga}{b}<\kappa_2, 0<d_2<d_{2*},\\
\sig_{1}(d_2 ,0,b\eta-\ga)~\text{changes sign at least once as}~d_1~\text{varies from}~0~\text{to}+\infty,~\text{if}~\kappa_1<\F{\ga}{b}<\kappa_2, d_2>d_{2*},\\
\sig_{1}(d_2 ,0,b\eta-\ga)<0,~\text{if}~\kappa_2<\F{\ga}{b}<\kappa_3,d_1\ge\hat{d_1},d_2>0,\\
\sig_{1}(d_2 ,0,b\eta-\ga)~\text{changes sign at least once as}~d_1~\text{varies from}~0~\text{to}+\infty,~\text{if}~\kappa_2<\F{\ga}{b}<\kappa_3,d_2<\bar{d_2},\\
\sig_{1}(d_2 ,0,b\eta-\ga)~\text{changes sign at least twice as}~d_1~\text{varies from}~0~\text{to}+\infty,~\text{if}~\kappa_2<\F{\ga}{b}<\kappa_3,d_2>\bar{d_2},\\
\sig_{1}(d_2 ,0,b\eta-\ga)<0,~\text{if}~\kappa_3<\F{\ga}{b}<\kappa_4,d_1\ge\hat{d_1},d_2>0,\\
\sig_{1}(d_2 ,0,b\eta-\ga)~\text{changes sign at least once as}~d_1~\text{varies from}~0~\text{to}+\infty,~\text{if}~\kappa_3<\F{\ga}{b}<\kappa_4,d_2\le\F{1}{\lam^{**}},\\
\sig_{1}(d_2 ,0,b\eta-\ga)~\text{changes sign at least twice as}~d_1~\text{varies from}~0~\text{to}+\infty,~\text{if}~\kappa_3<\F{\ga}{b}<\kappa_4,\F{1}{\lam^{**}}<d_2<\tilde{d_2},\\
\sig_{1}(d_2 ,0,b\eta-\ga)<0,~\text{if}~\kappa_3<\F{\ga}{b}<\kappa_4,d_2>\tilde{d_2},\\
\sig_{1}(d_2 ,0,b\eta-\ga)<0,~\text{if}~\F{\ga}{b}\ge\kappa_4.
\end{cases}
\eeq

\subsubsection{Existence of $q_0$}

Recall that if $\sig_1(d_2,0,b\eta-\ga)\le0$, then $\sig_1(d_2,q,b\th_1-\ga)<0$ for all $0<q<q^{*}(d_1,r)$;
if $\sig_1(d_2,0,b\eta-\ga)>0$,
then there exists a unique $\tilde{q}\in(0,q^{*}(d_1,r))$ such that $\sig_1(d_2,q,b\th_1-\ga)>0$ for $0<q<\tilde{q}$ and $\sig_1(d_2,q,b\th_1-\ga)<0$ for $\tilde{q}<q<q^{*}(d_1,r)$.
Hence we can define
\beq
\label{q_0}
q_0(d_1,d_2,\ga)=
\begin{cases}
0,~\text{if}~\sig_1(d_2,0,b\eta-\ga)\le0,\\
\tilde{q},~\text{if}~\sig_1(d_2,0,b\eta-\ga)>0.
\end{cases}
\eeq
Moreover, $\sig_1(d_2,q_0(d_1,d_2,\ga),b\th_1-\ga)=0$

\subsubsection{Properties of $q_0$ and proof of Theorem \ref{thm1}}

Since $\F{\pa\th_1}{\pa q}<0$,
it follows from Lemma \ref{lemma of eigenvalue}(iii)-(iv)
that $\F{\pa}{\pa\ga}q_0(d_1,d_2,\ga)<0$.
Define
\beq
\ga_0=\sup\limits_{\om\in H^{1}([0,1])\backslash\{0\}}
\F{-d_2\i_{0}^{1}\om^2_x\dx+b\i_0^1\th_1\om^2\dx}{\i_{0}^{1}\om^2\dx}.
\eeq
It is easy to see from the definition of $q_0$ that $q_0(d_1,d_2,\ga_0)=0$.
Recall from Lemma \ref{lemma of eigenvalue}(v)
and Lemma \ref{lemma of theta} that
\beq
\begin{cases}
\lim\limits_{d_1\rightarrow+\infty}\th_1=\F{\i_{0}^{1}r(x)\dx-q}{\i_{0}^{1}\f{r(x)}{K(x)}\dx},\\
\lim\limits_{d_2\rightarrow+\infty}\sig_1(d_2,q,b\th_1-\ga)=b\I_{0}^{1}\th_1\dx-q-\ga.
\end{cases}
\eeq
Then we obtain that
\beq
\label{g}
\begin{cases}
\lim\limits_{d_1\rightarrow+\infty}q_0(d_1,d_2,\ga)=\hat{q}_0,,\\
\lim\limits_{d_2\rightarrow+\infty}q_0(d_1,d_2,\ga)=\tilde{q}_0,\\
\lim\limits_{(d_1,d_2)\rightarrow(+\infty,+\infty)}q_0(d_1,d_2,\ga)=\F{b\i_0^1r(x)\dx-\ga\i_0^1\frac{r(x)}{K(x)}\dx}{b+\frac{r(x)}{K(x)}\dx},
\end{cases}
\eeq
where $\hat{q}_0$ and $\tilde{q}_0$ satisfy
\beq
\begin{cases}
\sig_1(d_2,\hat{q}_0,\F{b(\i_{0}^{1}r(x)\dx-\hat{q}_0)}{\i_{0}^{1}\f{r(x)}{K(x)}\dx}-\ga)=0,\\
b\I_0^1\th_1|_{q=\tilde{q}_0}\dx=\ga+\tilde{q}_0.
\end{cases}
\eeq
Now, we can complete the proof of Theorem \ref{thm1}
by means of \eqref{f}, \eqref{q_0} and \eqref{g}.

\subsection{Case $q_1\neq q_2$}
\label{case q1 neq q2}

In this subsection, we shall consider a general case where $q_1\neq q_2$,
which implies that the predators may experience different flow
speeds than the prey due to behavioral differences.
Similarly, we can see from \cite[Theorem 6.1(i)]{NWW} that $(0,0)$ is globally asymptotically stable
(among all nonnegative and nontrivial initial data)
provided that $q_1\ge q_1^*$,
where $q_1^*$ is defined in Lemma \ref{existence of theta}.
Hence we need  to investigate the sign of the principal eigenvalue $\sig_1(d_2,q_2,b\th_1-\ga)$.
Define
\beq
\label{mu1}
\mu_1(q_1,q_2):=\sig_1(d_2,q_2,b\th_1-\ga).
\eeq
Then it follows from Lemma \ref{lemma of theta}(iii) and \cite[Lemma 6.1]{NWW} that
\beq
\label{monotonicity of mu}
\begin{cases}
\mu_1(q_1,q_2)~\text{ is strictly decreasing with respect to}~q_1~\text{in}~[0,q_1^*),\\
\mu_1(q_1,q_2)~\text{ is strictly decreasing with respect to}~q_2~\text{in}~[0,+\infty),
\end{cases}
\eeq
provided that
$r'(x)\ge(\F{r(x)}{K(x)})'\max\limits_{x\in[0,1]}K(x)>0$ on $[0,1]$.

\subsubsection{Proof of Theorem \ref{thm2}}
\label{Proof of Theorem thm2}

Indeed, the sign of $\mu_1(0,0)$ has been considered in Section \ref{sign of mu(0,0)}.

{\ni\bf Proof of Theorem \ref{thm2}: }
Note that
\beq
\lim\limits_{q_1\rightarrow (q_1^*)^+}\mu_1(q_1,q_2)
=\sig_1(d_2,q_2,-\ga)=\sig_1(d_2,q_2,0)-\ga<0
\eeq
based on \cite[Proposition 2.1]{BDJS}.
On the other hand, it follows from Lemma \ref{properties of eta}(v) and
Lemma \ref{lemma of eigenvalue}(iii)(iv) that

\beq
\label{G}
\lim\limits_{q_2\rightarrow+\infty}\mu_1(0,q_2)=\lim\limits_{q_2\rightarrow+\infty}\sig_1(d_2,q_2,b\eta-\ga)<\lim\limits_{q_2\rightarrow+\infty}\sig_1(d_2,q_2,b\max\limits_{x\in[0,1]}K-\ga)<0.
\eeq

When $\mu_1(0,0)\le0$,
by the strict monotonicity of $\mu_1(q_1,q_2)$ with respect to $q_1,q_2$,
we conclude
that
\beq
\mu_1(q_1,q_2)<0~\text{for any}~(q_1,q_2)\in[0,q^{*}(d_1,r))\times[0,+\infty).
\eeq

When $\mu_1(0,0)>0$,
by the strict monotonicity of $\mu_1(q_1,q_2)$ again
and \eqref{G},
it indicates that
there exists a unique $\tilde{q}_2>0$ such that
\beq
\mu_1(0,q_2)
\begin{cases}
>0,~\text{for}~0<q_2<\tilde{q}_2,\\
=0,~\text{if}~q_2=\tilde{q}_2,\\
<0,~\text{for}~q_2>\tilde{q}_2.
\end{cases}
\eeq
For $q_2\ge\tilde{q}_2$,
by virtue of the strict monotonicity of $\mu_1(q_1,q_2)$ with respect to $q_1$ in $[0,q^{*}(d_1,r))$,
one can conclude that
\beq
\mu_1(q_1,q_2)<0~\text{for any}~(q_1,q_2)\in[0,q^{*}(d_1,r))\times[\tilde{q}_2,+\infty).
\eeq
While for $q_2<\tilde{q}_2$,
the strict monotonicity
of $\mu_1(q_1,q_2)$ with respect to $q_1$ yields that
there exists a unique continuous curve $q_1=\tilde{q}_1(q_2)$ defined in $q_2\in[0,\tilde{q}_2)$
such that
\beq
\mu_1(q_1,q_2)
\begin{cases}
>0,~\text{for}~0\le q_1<\tilde{q}_1(q_2), 0\le q_2<\tilde{q}_2,\\
=0,~\text{if}~q_1=\tilde{q}_1(q_2), 0\le q_2<\tilde{q}_2,\\
<0,~\text{for}~\tilde{q}_1(q_2)<q_1<q^{*}(d_1,r), 0\le q_2<\tilde{q}_2.
\end{cases}
\eeq
In particular,
$\tilde{q}_1(\tilde{q}_2)=0$,
$\F{\mathrm{d}\tilde{q}_1(q_2)}{\mathrm{d} q_2}<0$,
and $\mu_1(\tilde{q}_1(0),0)=0$.
Clearly, $\mu_1(\tilde{q}_1(0),0)=0$ implies that
$\tilde{q}_1(0)=q^{*}(d_1,r)$.
Hence, we obtain Theorem \ref{thm2}.
\proofend

\subsubsection{Proof of Theorem \ref{thm3}}
\label{Proof of Theorem thm3}

In this subsection,
we will consider the
influence of $b$ and the ratio $\tau=q_1:q_2$ of flow speeds experienced by predators and prey on the dynamical behavior of system \eqref{system}.
Now we set $q_1=\tau q_2$ in this subsection.
For simplicity,
we let $\sig_1(d_2,q_2,b\th_1-\ga)=\lam_1(\tau,q_2,b,\ga)$.

\pfthm3
In view of $q_1<q^{*}(d_1,r)$,
we arrive at $q_2<\F{q^{*}(d_1,r)}{\tau}$.
Since it follows from Lemma \ref{lemma of theta}(iii) that $(\th_1)_x>0$,
then
\beq
\begin{cases}
\lam_1(\tau,q_2,0,\ga)=\sig_1(d_2,q_2,-\ga)<0,\\
\lim\limits_{b\rightarrow+\infty}\lam_1(\tau,q_2,b,\ga)
\ge\lim\limits_{b\rightarrow+\infty}\sig_1(d_2,q_2,b\th_1(0)-\ga)=+\infty.
\end{cases}
\eeq
We hence get that there exists a unique $b_{\tau}(q_2)\in(0,+\infty)$
continuously depending on $q_2$
such that
\beq
\lam_1(\tau,q_2,b,\ga)
\begin{cases}
<0~\text{if}~b<b_{\tau}(q_2),\\
=0~\text{if}~b=b_{\tau}(q_2),\\
>0~\text{if}~b>b_{\tau}(q_2).
\end{cases}
\eeq
\proofend

In view of $\lam_1(\tau,q_2,b_{\tau}(q_2),\ga)=0$,
we have that there exist principal eigenfunctions $\phi,\varphi>0$ such that
\beq
\label{gg}
\left\{\arraycolsep=1.5pt
\begin{array}{lll}
d_2\varphi_{xx}-q_2\varphi_x+(b_{\tau}(q_2)\th_1-\ga)\varphi=0,&x\in(0,1),\\[2mm]
d_2\varphi_x(0)-q_2\varphi(0)=\varphi_x(1)=0,&
\end{array}
\right.
\eeq
and
\beq
\label{h}
\left\{\arraycolsep=1.5pt
\begin{array}{lll}
d_2\phi_{xx}+(b_{\tau}(0)\eta-\ga)\phi=0,&x\in(0,1),\\[2mm]
d_2\phi_x(0)-q_2\phi(0)=\phi_x(1)=0.&
\end{array}
\right.
\eeq
Then \eqref{h} implies that
\beq
b_{\tau}(0)=\inf\limits_{\om\in H^{1}([0,1])\backslash\{0\}}
\F{d_2\i_{0}^{1}\om^2_x\dx+\ga\i_0^1\om^2\dx}{\i_{0}^{1}\eta\om^2\dx}\ge\F{\ga}{\max\limits_{x\in[0,1]}K}.
\eeq
Set $w:=\F{\phi_x}{\phi}$.
Then by some straightforward computations,
one can find from Lemma \ref{int eta>kappa1} and \eqref{h} that
\beq
\left\{\arraycolsep=1.5pt
\begin{array}{lll}
-d_2w_{xx}-2d_2ww_x=b_{\tau}\varphi\eta_x>0,&x\in(0,1),\\[2mm]
w(0)=w(1)=0.&
\end{array}
\right.
\eeq
An application of the maximum principle yields
$w>0$.

Differentiating \eqref{gg}  with respect to $q_2$,
we have that
\beq
\label{i}
\left\{\arraycolsep=1.5pt
\begin{array}{lll}
d_2\varphi^{'}_{xx}-q_2\varphi^{'}_x-\varphi_x+(b_{\tau}(q_2)\th_1-\ga)\varphi^{'}+b^{'}_{\tau}(q_2)\th_1\varphi+b_{\tau}(q_2)\dot{\th_1}\tau\varphi=0,&x\in(0,1),\\[2mm]
d_2\varphi^{'}_x(0)-q_2\varphi^{'}(0)=\varphi(0),&\\[2mm]
\varphi^{'}_x(1)=0,
\end{array}
\right.
\eeq
where $\dot{\th_1}:=\F{\pa \th_1}{\pa(\tau q_2)}$.
Note that $\th_1|_{q_2=0}=\eta$,
$\varphi|_{q_2=0}=\phi$,
$\varphi_x|_{q_2=0}=\phi_x>0$.
Multiplying \eqref{gg} by $\varphi^{'}e^{-\frac{q_2}{d_2}x}$ and \eqref{i} by $\varphi e^{-\frac{q_2}{d_2}x}$, subtracting the resulting equations and then integrating over $(0,1)$,
one finally derives that
\beq
b^{'}_{\tau}(q_2)\I_0^1\th_1\varphi^2e^{-\frac{q_2}{d_2}x}\dx
=-b_{\tau}(q_2)\tau\I_0^1\varphi^2e^{-\frac{q_2}{d_2}x}\dot{\th_1}\dx
+\F{q_2}{2d_2}\I_0^1\varphi^2e^{-\frac{q_2}{d_2}x}\dx
+\F12\varphi(1)e^{-\frac{q_2}{d_2}}+\F12\varphi(0)>0.
\eeq
Let $\xi$ be the principal eigenfunction of the following problem
\beq
\label{xi}
\left\{\arraycolsep=1.5pt
\begin{array}{lll}
d_2\xi_{xx}-q_2\xi_x+\ga\xi=\Gamma_1\xi,&x\in(0,1),\\[2mm]
d_2\xi_x(0)-q_2\xi(0)=\xi_x(1)=0.&
\end{array}
\right.
\eeq
Clearly, $\Gamma_1<\ga$.
Multiplying \eqref{gg} by $\xi e^{-\frac{q_2}{d_2}x}$ and \eqref{xi} by $\varphi e^{-\frac{q_2}{d_2}x}$, subtracting the resulting equations and then integrating over $(0,1)$,
one can obtain that
\beq
b_{\tau}(q_2)=
\F{(2\ga-\Gamma_1)\i_0^1\varphi\xi  e^{-\frac{q_2}{d_2}x}\dx}{\i_0^1\varphi\xi  e^{-\frac{q_2}{d_2}x}\th_1\dx}
\ge\F{(2\ga-\Gamma_1)}{\max_{x\in[0,1]}\th_1}.
\eeq
This indicates from $\lim\limits_{q_2\rightarrow(\frac{q^{*}(d_1,r)}{\tau})^{-}}\th_1=0$ that
\beq
\lim\limits_{q_2\rightarrow(\frac{q^{*}(d_1,r)}{\tau})^{-}}b_{\tau}(q_2)=+\infty.
\eeq
Therefore, we finish the proof of Theorem \ref{thm3}.
\proofend

%
%



%

\end{document}